\documentclass[12pt,letterpaper]{article}
\usepackage{amsfonts,amssymb,latexsym}

\newtheorem{theorem}{Theorem}[section]
\newtheorem{lemma}{Lemma}[section]

\newtheorem{definition}{Definition}[section]

\newcommand{\eop}{$\quad\spadesuit$}
\newenvironment{proof}{{\bf Proof.}}{\eop\medskip}

\begin{document}

\author{Ramiro de la Vega}
\title{Homogeneity properties in large compact $S$-spaces}
\maketitle

\begin{abstract}
Under $\Diamond $ we show how to construct a large compact $S$-space while having some control over its group of autohomeomorphisms. In particular we can make the space rigid or $h$-homogeneous (i.e. any two clopen subsets are homeomorphic). We also get a space which is $B$-homogeneous and not $h$-homogeneous. A topological space is $B$-homogeneous if it has a base every element of which can be mapped onto every other one by a homeomorphism of the whole space.
\end{abstract}

\section{Introduction}

A space is {\it hereditarily separable} ($HS$) if every subspace is separable. An $S${\it-space} is a regular $HS$ space with a non-Lindel\"of subspace. A space $X$ is {\it homogeneous} if for every $x,y\in X$ there is a homeomorphism $f$ of $X$ onto $X$ with $f\left( x\right) =y$. In \cite{fed}, Fedorchuk constructed under $\Diamond $ a compact $S$-space of size $2^{\aleph _1}$ (in fact with every point of the space having character $\aleph _1$). There are compact homogeneous $S$-spaces under $CH$ (see \cite{yo}) but by classic results such a space must have size at most $2^{\aleph_0}$. In this paper we show how to construct a Fedurchuk space (i.e. compact, HS and with no non-trivial convergent sequences) keeping control over its group of autohomeomorphisms. This allows us to get such spaces with various homogeneity properties.

A topological space is $B${\it-homogeneous} if it has a base every element of which can be mapped onto every other one by a homeomorphism of the whole space. A zero-dimensional space is $h${\it -homogeneous} if any two clopen subsets are homeomorphic to each other. It is shown in \cite{mat} that any $B$-homogeneous zero-dimensional first-countable compact space is $h$-homogeneous (and hence homogeneous). Matveev also constructed an example of a pseudocompact zero-dimensional $B$-homogeneous space which is not $h$-homogeneous, and he asked whether any $B$-homogeneous zero-dimensional compact space is $h$-homogeneous. We show that under $\Diamond $ the answer is no. In fact our example is a Fedurchuck space. We can also get an $h$-homogeneous Fedurchuck space and a rigid one (i.e. one for which the only autohomeomorphism of the space is the identity map).

We use an inverse limit approach similar to the one used in \cite{mirna}, where $\Diamond $ is used to capture (and eventually destroy) all potential left separated $\omega _1$-sequences, making the space $HS$. In order to keep control over the autohomeomorphisms of the space we need to do two things: we need to destroy some potential autohomeomorphisms (for which we use $\Diamond $ again), but we also want to preserve some of them, for which we use a ``lifting" technique used in \cite{kun} for getting certain kind of $L$-spaces. In Section \ref{terms} we introduce some terminology and give basic facts used during the construction. In Section \ref{lifting} we introduce the lifting technique mentioned before, which will be used at the succesor stages of the inverse limit. The main construction is done in Section \ref{main}; we recommend the reader to go straight to that section and go back to previous sections as needed. Finally in Section \ref{apps} we explain briefly how to use the main construction to get the various spaces
we mentioned at the beginning.

\section{General Facts and Terminology\label{terms}}

Suppose $G$ is a group acting on a topological space $X$.

\begin{definition}
We say that $f\in Hom(X)$ is {\it locally in} $G$ if for every $x\in X$ there is an open neighborhood $U$ of $x$ and a $g\in G$ such that $f(u)=g(u)$ for all $u\in U$. We write $H(X,G)$ for the set of autohomeomorphisms of $X$ which are locally in $G$.
\end{definition}

Note that if $X$ is compact and zero-dimensional then $f$ is locally in $G$ if and only if there is a finite partition $\{P_i:i\in I\}$ of $X$ into clopen subsets and there is a finite $\{g_i:i\in I\}\subseteq G$ such that $f\restriction P_i=g_i$ for all $i\in I$. Also note that $H(X,G)$ is actually a subgroup of $Hom(X)$ and that $H(X,H(X,G))=H(X,G)$.

\begin{definition}
Two subsets $A,B\subseteq X$ are $G$-{\it equivalent} if there is an $f\in H(X,G)$ such that $fA=B$.
\end{definition}

\begin{lemma}\label{base}
Suppose that $U\subseteq X$ is clopen and let ${\cal B}$ be the set of all clopen subsets of $X$ which are $G$-equivalent to $U$. If ${\cal B}$ is a $\pi$-base for $X$ and $X=\bigcup {\cal B}$ then ${\cal B}$ is a base for $X$.
\end{lemma}

\begin{proof}
Fix a clopen $V\subseteq X$ and $x\in V$. We want to find an element of ${\cal B}$ contained in $V$ and containing $x$. Without loss of generality there is a $B\in {\cal B}$ such that $B\subseteq X\setminus V$. Since $X=\bigcup {\cal B}$, there is an $f\in H(X,G)$ such that $x\in fB$, and we can fix $W\subseteq B$ clopen such that $x\in f(W) \subset V$. Since $B$ is a $\pi$-base, there is an $f'\in H(X,G)$ such that $f'(B) \subseteq V\setminus f(W)$. Now we define $h\in H(X,G)$ as follows: if $x\in W$ we let $h(x) =f(x)$, if $x\in B\setminus W$ we let $h(x)=f'(x)$, if $x\in f(W)$ we let $h(x)=f^{-1}(x)$, if $x\in f'(B\setminus W)$ we let $h(x)=f'^{-1}(x)$ and we let $h(x)=x$ otherwise. Then clearly $h\in H(X,G)$ and $x\in hB\subseteq V$. Thus $hB$
is the desired element of ${\cal B}$.
\end{proof}

\begin{definition}
Given a countable sequence $S$ of closed subsets of the Cantor space $C$ and given closed $A_0,A_1\subseteq C$, we say that $(A_0,A_1)$ {\it preserves} $S$ if whenever $x\in A_0\cap A_1$ is a strong limit point of $S$ and $U$ is an open neighborhood of $x$, there are $n_0,n_1\in dom(S)$ with $S(n_0) \subseteq U\setminus A_1$ and $S(n_1) \subseteq U\setminus A_0$.
\end{definition}

\begin{definition}
Given a countable set ${\cal S}$ of countable sequences of closed subsets of the Cantor space $C$ and given a convergent sequence $B=\{b_n:n\in\omega\} \subseteq C$ with limit point $b$; we say that $B$ {\it avoids} ${\cal S}$, if whenever $b$ is a strong limit point of $S\in {\cal S}$ and $U$ is an open neighborhood of $b$, there is an $n\in dom(S)$ such that $S(n) \subseteq U\setminus (B\cup\{b\})$.
\end{definition}

Note in particular that any constant sequence avoids any set ${\cal S}$ and hence satisfies the hypothesis of the next lemma.

\begin{lemma}\label{split}
Suppose ${\cal S}$ is a countable set of countable sequences of closed subsets of the Cantor space $C$. Suppose also that $B=\{b_n:n\in \omega \} \subseteq C$ is a convergent sequence which avoids ${\cal S}$ and let $b$ be its limit point. Then there are closed $A_0,A_1\subseteq C$ such that i) $A_0\cup A_1=C$, ii) $A_0\cap A_1=\{b\}$, iii) $B\subseteq A_0$ and iv) $(A_0,A_1)$ preserves every sequence in ${\cal S}$.
\end{lemma}

\begin{proof}
This result was essentially proved in \cite[proof of Theorem 0.1]{kun}. The difference is that there they don't have the sequence $B$ to worry about. Thus can assume here that the $b_n$'s are all distinct and all different from $b$. Note that the order type of the sequences in ${\cal S}$ is irrelevant so we may assume that ${\cal S}=\{\langle S_n^i:n\in\omega \rangle :i\in \omega \}$. Fix a decreasing sequence $\langle V_m:m\in \omega \rangle$ of clopen subsets of $C$ such that $V_0=C$ and $\bigcap_{m\in \omega}V_m=B\cup \{b\}$. We will construct inductively a decreasing sequence $\langle W_m:m\in\omega\rangle$ of clopen subsets of $C$ such that $W_0=C$, $W_m\subseteq V_m$ and $b_m\in W_{2m}\setminus W_{2m+1}$ for all $m\in\omega$. Note that $\bigcap_{m\in \omega}W_m=\{b\}$ and that if we let $A_0=\{b\} \cup \bigcup_{m\in \omega}(W_{2m}\setminus W_{2m+1})$ and $A_1=\{b\} \cup\bigcup_{m\in \omega}(V_{2m+1}\setminus V_{2m+2})$, then i), ii) and iii) are satisfied. Now we explain how to construct the sequence so that iv) is also satisfied.

Fix a local base $\{O_j:j\in \omega\}$ for $b$ in $C$ and a function $(\psi,\mu):\omega\rightarrow \omega \times \omega$ such that $(\psi,\mu)^{-1}(i,j)$ contains at least one odd and one even natural number for each $(i,j)\in\omega \times \omega$. If $W_{2m}$ is already defined and $b$ is a strong limit point of $\langle S_n^{\psi (2m)}:n\in \omega\rangle$ then (since $B$ avoids ${\cal S}$) there is an $n\in \omega$ such that $S_n^{\psi (2m)}\subseteq W_{2m}\cap O_{\mu(2m)}$ and $S_n^{\psi (2m)}\cap (B\cup \{b\})=\emptyset$. Thus there is a $k\in \omega$ and a clopen neighborhood $U$ of $b_m$ such that $S_n^{\psi (2m)}\subseteq W_{2m}\setminus (V_k\cup U)$ and $U\cap B=\{b_m\}$; then we let $W_{2m+1}=(W_{2m}\cap V_k) \setminus U$. If $b$ wasn't a strong limit point of $\langle S_n^{\psi (2m)}:n\in \omega \rangle$ then we just let $W_{2m+1}=(W_{2m}\cap V_{2m+1}) \setminus U$ where $U$ is any clopen with $U\cap B=\{b_m\}$. Note that in either case we have that $\{b_k:k>m\} \subseteq W_{2m+1}\subseteq V_{2m+1}\cap W_{2m}$ and $b_m\notin W_{2m+1}$. If $W_{2m+1}$ is already defined we follow the same steps (now working with $S^{\psi (2m+1)}$ and $O_{\mu(2m+1)}$) except that now we want to keep $b_{m+1}$ inside $W_{2m+2}$, so we just let $W_{2m+2}=W_{2m+1}\cap V_k$ in the first case and $W_{2m+2}=W_{2m+1}\cap V_{2m+2}$ in the second case.

Now suppose that $b$ is a strong limit point of $\langle S_n^i:n\in\omega \rangle$ and $O$ is an open neighborhood of $b$. Then there is a $j\in \omega $ such that $b\in O_j\subseteq O$ and there are $n_0$ and $n_1$ such that $(\psi ,\mu )(2n_0) =( \psi ,\mu )(2n_1+1)=(i,j)$. It is clear from the construction that these $n_0$ and $n_1$ are as required in the definition of preserving $\langle S_n^i:n\in \omega \rangle $.
\end{proof}

We close this section with a result which was essentially proven in \cite[Lemma 6.4]{kun}. There instead of two clopen sets $K$ and $L$ they had points $a$ and $b$ of the Cantor space and the requirement at the end was $ga=b$. However an almost identical argument shows the following.

\begin{lemma}\label{free}
Suppose $G$ is a countable subgroup of autohomeomorphisms of the Cantor space $C$ acting freely on $C$ and suppose $K$ and $L$ are disjoint non-empty clopen subsets of $C$. Then there is an autohomeomorhism $g$ of $C$ such that $gK=L$ and $\langle G,g\rangle$ acts freely on $C$.
\end{lemma}

\section{Lifting Homeomorphisms}\label{lifting}

Suppose $G$ is a countable group acting freely on the Cantor space $C=2^\omega $. Given two closed $A_0,A_1\subseteq C$ and $p\in C$ such that $A_0\cap A_1=\{p\}$ and $A_0\cup A_1=C$, we let $Y=Y(A_0,A_1,G)$ be the space 
$$
Y=\left\{(x,\varphi)\in C\times 2^G:(\forall g\in G)\left( x\in gA_{\varphi (g)}\right)\right\}.
$$
We let $G$ act on $2^G$ by $(g\varphi)(h)=\varphi (g^{-1}h)$ where $g$ and $h$ denote elements
of $G$ while $\varphi$ and $g\varphi$ are in $2^G$. Furthermore $G$ acts on $Y$ by $g(x,\varphi)=(gx,g\varphi)$.

For any choice of $A_i$'s we have that $Y$ is a closed subspace of $C\times 2^G$ and the action described above is free and continuous. We also have that the projection into the first coordinate $\pi:Y\rightarrow C$ is surjective.

\begin{lemma}
Suppose that ${\cal P}$ is a finite partition of $Y$ into clopen subsets. Then there is a finite partition ${\cal Q}$ of $C$ into clopen subsets and a function $\Psi:{\cal Q}\rightarrow [{\cal P}]^{\leq 2}$ such that i) $\pi^{-1}Q\subseteq \bigcup \Psi (Q)$ for all $Q\in {\cal Q}$ and ii) whenever $|\Psi(Q)| =2$, there is a unique $x\in Q$ such that $\Psi(Q)$ splits $\pi^{-1}(x)$.
\end{lemma}

\begin{proof}
For each clopen $Q\subseteq C$ let $\Psi (Q)=\{P\in{\cal P}:P\cap\pi^{-1}Q\neq\emptyset\}$. Since the action of $G$ is free we have that $|\pi^{-1}(x)|\leq 2$ for all $x\in C$. Thus given $x\in C$ there is a clopen neighborhood $Q_x$ of $x$ such that $|\Psi (Q)|\leq2$. Now using compactness of $C$ and taking boolean combinations we can get a finite partition ${\cal Q}'$ of $C$ with $|\Psi(Q)|\leq2$ for all $Q\in{\cal Q}'$. Now suppose that $Q\in {\cal Q}'$ is such that $|\Psi (Q)| =2$. Since the two clopens in $\Psi (Q)$ only mention finitely many elements of $G$ and $\Psi (Q)$ splits $\pi^{-1}(x)$ only if (maybe not even if) $x=gp$ for some $g$ mentioned, we can further partition $Q$ to isolate the finitely many such possible $x$'s. We then let ${\cal Q}$ to be the corresponding refinement of ${\cal Q}'$.
\end{proof}

Note that by the same argument, we can in fact get ${\cal Q}$ in the previous lemma so that it refines any finite partition of $C$ given in advance.

The next result shows that non-$G$-equivalent clopens of $C$ remain non-$G$-equivalent in $Y$.

\begin{lemma}\label{nonequipres}
Suppose that $U,V$ are disjoint clopen subsets of $C$ such that $\pi ^{-1}U$ is $G$-equivalent to $\pi ^{-1}V$. Then $U$ and $V$ are $G$-equivalent.
\end{lemma}

\begin{proof}
Fix $f\in H(Y,G)$ such that $f(\pi^{-1}U)=\pi^{-1}V$. Without loss of generality we may assume that $f^2=id$ and that $f\restriction Y\setminus \pi^{-1}(U\cup V)=id$. Fix a finite partition ${\cal P}$ of $Y$ into clopen sets so that $f$ acts as a different element of $G$ in each element of the partition. Let ${\cal Q}$ be a partition of $C$ refining $\{U,V,C\setminus (U\cup V)\}$ as in the previous lemma with the corresponding $\Psi :{\cal Q}\rightarrow [{\cal P}] ^{\leq 2}$.

Let $S=\{ x\in C:\pi^{-1}(x)$ is split by ${\cal P}\}$. Since $|S\cap Q|\leq 1$ for each $Q\in{\cal Q}$, we have that $S$ is finite. On the other hand, since the action of $G$ in $Y$ commutes with $\pi$, it is easy to see that $f(\pi^{-1}S)=\pi^{-1}S$ so we get that $|S\cap V|=|S\cap U|$. Then we can write $S\cap U=\{g_ip:i\in n\}$ and $S\cap V=\{h_jp:j\in n\}$ for some $n\in \omega$ and some $g_i,h_j$'s in $G$.

Now fix $W$ a clopen neighborhood of $p$ such that each $g_iW$ and each $h_jW$ is contained in an element of ${\cal Q}$. Let $W_U=\bigcup_{i\in n}g_iW$ and $W_V=\bigcup_{i\in n}h_iW$. We claim that $f(\pi^{-1}W_U) =\pi^{-1}W_V$.

Suppose that $(x,\varphi) \in \pi^{-1}W_U$ and fix $i\in n$ and $Q\in {\cal Q}$ such that $x\in g_iW\subseteq Q$. Let $\psi \in 2^G$ such that $( g_ip,\psi)$ and $(x,\varphi)$ belong to the same element $P$ of $\Psi(Q)$ and let $g\in G$ such that $f\restriction P=g\restriction P$. Now fix $j\in n$ such that $\pi (f(g_ip,\psi)) =h_jp$. Thus $gg_ip=h_jp$ and since the action of $G$ is free we get that $gg_i=h_j$ and therefore $\pi (f(x,\varphi))=\pi (gx,g\varphi)=gx=h_jg_i^{-1}x\in h_jW$. This shows that $f(x,\varphi) \in \pi^{-1}W_V$, hence proving that $f(\pi ^{-1}W_U) \subseteq \pi^{-1}W_V$. By a similar argument it follows that $f(\pi^{-1}W_V) \subseteq \pi ^{-1}W_U$ and since $f^2=id$ we get $f(\pi ^{-1}W_U) =\pi^{-1}W_V$.

Finally we apply the last lemma one more time to get a partition ${\cal R}$ of $C\setminus (W_U\cup W_V)$ such that for each $R\in {\cal R}$ there is a $g_R\in G$ with $f\restriction \pi^{-1}R=g_R\restriction \pi^{-1}R$. Now we can define $h:C\rightarrow C$ as follows: if $x\in R\in {\cal R}$ let $h(x) =g_R(x)$, if $x\in g_iW$ let $h(x)=h_ig_i^{-1}(x)$ and if $x\in h_jW$ let $h(x)=g_jh_j^{-1}(x)$. It is clear that $h$ is an autohomeomorphism of $C$, $h$ is locally in $G$ and $hU=V$.
\end{proof}

\begin{definition}
Suppose $(x,\varphi) \in Y$ and $|\pi^{-1}(x)|=2$. Then $x=gp$ for some (unique) $g\in G$. We say that $(x,\varphi)$ is of $0$-{\it type} if $\varphi (g)=0$ and $(x,\varphi)$ is of $1$-{\it type} if $\varphi (g)=1$.
\end{definition}

The next result shows that not all the clopens of $Y$ are $G$-equivalent.

\begin{lemma}\label{nonequiexists}
Let $K=\{(x,\varphi) \in Y:\varphi(e) =0\}$ where $e\in G$ is the identity element. Then $K$ and $
Y\setminus K$ are not $G$-equivalent.
\end{lemma}

\begin{proof}
By contradiction, fix a locally in $G$ autohomeomorphism $f$ of $Y$ such that $fK=Y\setminus K$. Without loss of generality we may assume that $f^2=id$. Let $S$ be the set of all $x\in C$ such that $\pi ^{-1}(x)$ has two points and $f$ is defined by different elements of $G$ at those two points. Then $f(\pi^{-1}S)=\pi^{-1}S$ and in fact $f(\pi^{-1}S\cap K)=\pi^{-1}S\cap (Y\setminus K)$. Note that $p\in S$ and that for any $x\in S$ with $x\neq p$ we either have $\pi ^{-1}(x) \subseteq \pi^{-1}S\cap K$ or $\pi^{-1}(x) \subseteq \pi^{-1}S\cap (Y\setminus K)$. Thus $\pi^{-1}S\cap K$ has $n+1$ points of $0$-type and $n$ points of $1$-type while $\pi^{-1}S\cap (Y\setminus K)$ has $n$ points of $0$-type and $n+1$ points of $1$-type, where $|S|=2n+1$. But this is impossible since the action of $G$ preserves the type of the points.
\end{proof}

\begin{lemma}\label{preserve}
Suppose that $S$ is a sequence of closed subsets of $C$ and $(A_0,A_1)$ preserves $gS$ for each $g\in G$. Suppose also that $x\in C$ is a strong limit point of $S$. Then every point in $\pi^{-1}(x)$ is a strong limit point of $\pi^{-1}S$.
\end{lemma}

\begin{proof}
Fix $(a,\varphi) \in Y$ and suppose that $a$ is a strong limit point of $S=\langle S_i:i\in I\rangle$. Fix $U\times W$ an open neighborhood of $(a,\varphi)$ which may be assumed to be of the form $W=\{\psi \in 2^G:\sigma \subseteq \psi\}$ for some $\sigma \subseteq \varphi$ finite. If $g\in dom(\sigma)$ and $a\neq gp$ then $a\notin gA_{1-\sigma(g)}$ so there is an open $V_g$ with $a\in V_g\subseteq U$ and $V_g\cap gA_{1-\sigma(g)}=\emptyset$. Let $V=\bigcap \{V_g:g\in dom(\sigma)$ and $a\neq gp\}$. Now we consider two cases:

i) If $a=gp$ for some $g\in dom(\sigma)$ (since $G$ acts freely, there is at most one such $g$) then $p=g^{-1}a\in g^{-1}V$ and $p$ is a strong limit point of $g^{-1}S$. Therefore (since $(A_0,A_1)$ preserves $g^{-1}S$) there is an $i\in I$ such that $g^{-1}S_i\subseteq g^{-1}V\setminus A_{1-\sigma(g)}$. Thus $S_i\subseteq V\setminus gA_{1-\sigma(g)}$.

ii) if $a\neq gp$ for all $g\in dom(\sigma)$ then we just use the fact that $a$ is a strong limit point of $S$ to get an $i\in I$ for which $S_i\subseteq V\setminus \{a\}$.

In either case we get that $\pi ^{-1}S_i\subseteq (U\times W)\setminus \{(a,\varphi)\}$, showing that $(a,\varphi)$ is a strong limit point of $\pi ^{-1}S$.
\end{proof}

\section{Main Construction}\label{main}

We start with a fixed countable group $G_0$ of autohomeomorphisms of the Cantor space $C=2^\omega $ which acts freely (i.e. $g(x)\neq x$ for all $x\in 2^\omega $ and for all $g\in G_0\setminus \{id\}$). We construct a space $X\subseteq C\times C^{\omega _1}$, a group $G\supseteq G_0$ and a continuous action $*:G\times X\rightarrow X$ such that (here $\pi :C\times C^{\omega _1}\rightarrow C$ is the projection onto the first coordinate):

i) $X$ is compact, has no isolated points and $\pi\restriction X$ is surjective.

ii) $X$ is hereditarily separable and has no non-trivial convergent
sequences.

iii) The action $*$ is free and for all $x\in X$ and $g\in G_0$, $g(\pi(x)) =\pi(g*x)$.

iv) $Hom(X) \subseteq H(X,G)$

v) $\{\pi^{-1}K:K\subseteq 2^\omega$ is clopen $\}$ forms a $\pi$-base of $X$.

We will define inductively $X_\alpha \subseteq C\times C^\alpha$ for $\alpha \leq \omega _1$ and then let $X=X_{\omega _1}$. We also construct groups $G_\alpha $ and actions $*_\alpha:G_\alpha \times X_\alpha \rightarrow X_\alpha $ and then let $G=G_{\omega _1}$ and $*=*_{\omega _1}$. Actually, here we only give a guideline for constructing the $G_\alpha$'s, leaving enough freedom so that we can use this same construction for the various applications.

For $\alpha \leq \beta \leq \omega _1$ let $\pi_\alpha ^\beta :C\times C^\beta \rightarrow C\times C^\alpha $ be the obvious projection (we will often drop the superscript if $\beta =\omega _1$). To start just let $X_0=C\times 1$ and $*_0$ is the natural action of $G_0$ induced in $X_0$. We want the following conditions to hold on each step of the induction:\\

{\bf C1)} $X_\alpha $ is closed in $C\times C^\alpha $ and has no isolated points.

{\bf C2)} $\pi_\beta^\alpha (X_\alpha) =X_\beta$ for every $\beta \leq \alpha $.

{\bf C3)} Each $*_\alpha $ is a free continuous action of $G_\alpha $ on $X_\alpha $.

{\bf C4)} If $\beta \leq \alpha $ then $G_\beta \subseteq G_\alpha $ and $\pi _\beta ^\alpha (g*_\alpha x) =g*_\beta ( \pi _\beta^\alpha (x))$ for all $x\in X_\alpha$ and $g\in G_\beta$.\\

This conditions almost determine the construction at limit ordinals. More precisely if $\gamma $ is a limit ordinal then $X_\gamma=\bigcap_{\alpha \in \gamma}(\pi _\alpha ^\gamma)^{-1}(X_\alpha)$, we let $G_\gamma '=\bigcup_{\alpha \in \gamma }G_\alpha$ and for any $g\in G_\gamma '$ and $x\in X_\gamma$, $g*_\gamma x$ is determined by the equations $\pi _\alpha ^\gamma (g*_\gamma x)=g*_\alpha (\pi_\alpha^\gamma (x))$ for $\alpha _0\leq \alpha<\gamma $ (where $\alpha _0=\min \{\alpha :g\in G_\alpha \}$). Then for each particular application, we will choose some $G_\gamma \supseteq G_\gamma '$ and extend the action $*_\gamma $ to $G_\gamma $. Note then that i) and iii) will be satisfied. To make sure that ii) and iv) are also satisfied we need $\Diamond$.

Assuming $\Diamond$ we can fix $\{ z_\gamma^\alpha :\alpha \in \omega_1 ,\gamma \in \alpha \}$ such that $\langle z_\gamma^\alpha :\gamma \in \alpha \rangle$ is an $\alpha $-sequence in $C\times C^\alpha$ for each $\alpha \in \omega _1$ and for every $\omega _1$-sequence $\langle x_\gamma :\gamma \in \omega _1\rangle$ in $C\times C^{\omega _1}$, the set $\{\alpha \in \omega _1:(\forall\gamma \in \alpha) \pi _\alpha (x_\gamma)=z_\gamma^\alpha \} $ is stationary. Now we can state the following condition which will ensure that $X$ is hereditarily separable (see \cite{mirna} for more details):\\

{\bf C5)} For all $\beta \leq \alpha \leq \omega _1$, if the $z_\gamma^\beta $ 's are all distinct and in $X_\beta $ then for all $x\in X_\alpha $ for which $\pi _\beta ^\alpha (x)$ is a limit point of $\langle z_\gamma ^\beta :\gamma \in \beta \rangle $ we have that $x$ is a strong limit point of $\langle X_\alpha \cap (\pi _\beta^\alpha)^{-1}(z_\gamma ^\beta) :\gamma \in \beta
\rangle $.\\

We want $X$ to have no convergent sequences. For this we can use our $\Diamond $-sequence to state the condition:\\

{\bf C6)} If $\{ z_n^\alpha :n\in \omega \}$ is a non-trivial convergent sequence in $X_\alpha$ with limit point $x\in X_\alpha $ then $|X_{\alpha +1}\cap (\pi _\alpha^{\alpha +1})^{-1}(
x)| >1$.\\

If $\langle x_n:n\in \omega \rangle $ was a convergent sequence in $X$ then we can extended it to a sequence $\langle x_\alpha :\alpha\in \omega _1\rangle $ in such a way that $\langle x_n:n\in\omega \rangle$ and $\langle x_\alpha :\alpha \geq \omega\rangle$ have disjoint clousures. Then for some $\alpha $ we would have that $\langle \pi _\alpha (x_n) :n\in \omega\rangle=\langle z_n^\alpha :n\in \omega \rangle $ is a non-trivial convergence sequence with limit point say $x$. By C5) we have that any element of $(\pi _\alpha)^{-1}(x)$ is a limit point of $\langle x_n:n\in\omega \rangle $. On the other hand by C6) we have that $|(\pi _\alpha)^{-1}(x)|>1$, which is a contradiction.

If $h$ is a homeomorphism from $X$ onto $X$, there is a club $Q_h\subseteq\omega _1$ and for each $\alpha \in Q_h$ there is a homeomorphism $h_\alpha$ from $X_\alpha $ onto $X_\alpha $ such that $\pi _\alpha (h(x))=h_\alpha (\pi_\alpha(x))$ for all $x\in X$. On the other hand any continuous function is determined by its restriction to a dense subset. Thus we can use our $\Diamond$-sequence to capture continuous functions. Formally we make the requirement:\\

{\bf C7)} If $\{z_n^\alpha :n\in \omega \}$ is a dense subset of $X_\alpha $, $f$ is a homeomorphism from $X_\alpha$ onto $X_\alpha$ such that $f(z_n^\alpha)=z_{\omega +n}^\alpha$ for $n\in\omega $ (note that there is at most one such $f$) and $f\notin H(X_\alpha ,G_\alpha)$, then there is no autohomeomorphism $h$ of $X$ such that $\alpha \in Q_h$ and $h_\alpha =f$.\\

To satisfy C7) we construct for each $\alpha \in \omega _1$ (for which the hypothesis of C7) hold) two sequences $\langle a_n^\alpha :n\in\omega \rangle$ and $\langle b_n^\alpha :n\in \omega
\rangle $ in $X_\alpha $ converging to $a^\alpha $ and $b^\alpha$ respectively, satisfying the following requirements:\\

{\bf R1)} There is an $x\in X_{\alpha +1}\cap ( \pi _\alpha ^{\alpha+1}) ^{-1}(b^\alpha)$ which is not a strong limit point of the sequence $\langle X_{\alpha +1}\cap ( \pi _\alpha ^{\alpha
+1}) ^{-1}(b_n^\alpha) :n\in \omega \rangle$.

{\bf R2)} For all $\beta \leq \alpha \leq \omega _1$ and for all $x\in X_\alpha \cap (\pi _\beta ^\alpha)^{-1}(a^\beta)$ we have that $x$ is a strong limit point of $\langle X_\alpha \cap (\pi_\beta ^\alpha)^{-1}(a_n^\beta) :n\in \omega\rangle$.

{\bf R3)} $f(a_n^\alpha)=b_n^\alpha$ for all $n\in \omega $.\\

Note that if we can meet the last three requirements then C7) is automatically satisfied because if $h$ was a counterexample, it would have to send the sequence $\langle X\cap (\pi _\alpha)^{-1}(a_n^\alpha):n\in \omega \rangle$ into the sequence $\langle X\cap (\pi _\alpha)^{-1}(b_n^\alpha):n\in \omega \rangle$ and the set $X\cap (\pi _\alpha)^{-1}(a^\alpha)$ into the set $X\cap (\pi_\alpha)^{-1}(b^\alpha)$. But this is impossible because any point in $X\cap (\pi_\alpha)^{-1}(a^\alpha)$ is a strong limit point of the first sequence while there are points in $X\cap (\pi_\alpha)^{-1}(b^\alpha)$ which are not a strong limit points of the second sequence.

Finally, to make sure that v) holds at the end (and in fact at each stage of the construction), we fix in advance an enumeration $\{ W_n:n\in \omega\}$ of the clopen subsets of $C\times 1$ and make the requirement:\\

{\bf C8)} For each $\alpha \leq \omega _1$ and for all $x\in X_\alpha $, $x$ is a strong limit point of $\langle X_\alpha \cap (\pi_0^\alpha)^{-1}(W_n):n\in \omega \rangle$.\\

All conditions are automatically satisfied at limit stages.

For the succesor stage, fix $\alpha \in \omega _1$ and assume $X_\beta $, $G_\beta $ and $*_\beta $ have been constructed for all $\beta \leq \alpha $. Now we define $X_{\alpha +1}$ and we show how to lift the action of $G_\alpha$($=G_{\alpha +1}'$) to $X_{\alpha +1}$. Then in each particular application of the construction we will define $G_{\alpha +1}\supseteq G_{\alpha+1}' $ and extend the action.

Let ${\cal S}''$ be the set of all sequences of the form $\langle g(X_\alpha \cap (\pi_\beta^\alpha)^{-1}(z_\gamma^\beta)):\gamma \in \beta \rangle$ with $g\in G_\alpha $ and $\beta \leq\alpha $ such that the $z_\gamma^\beta$ 's are all distinct and in $X_\beta$, together with all the sequences of the form $\langle g(X_\alpha \cap (\pi_\beta^\alpha)^{-1}(a_n^\beta)):n\in \omega
\rangle$ with $g\in G_\alpha$ and $\beta <\alpha$ ($\beta$ for which condition C7) was considered), together with the sequence $\langle X_\alpha \cap (\pi _0^\alpha)^{-1}(W_n):n\in \omega \rangle$.

If $f:X_\alpha \rightarrow X_\alpha $ satisfies the hypothesis of C7) fix $a^\alpha \in X_\alpha $ witness of the fact that $f\notin H(X_\alpha,G_\alpha)$ (i.e. there is no clopen neighborhood of $a^\alpha$ on which $f$ acts as an element of $G_\alpha $) and let $b^\alpha =f(a^\alpha)$. If there is a $g_0\in G_\alpha $ such that $g_0*_\alpha a^\alpha =b^\alpha $ then find a sequence $\langle O_n:n\in \omega\rangle$ of clopen subsets of $X_\alpha$ converging to $a^\alpha$ such that for all $n\in \omega $ and all $x\in O_n$ we have that $g_0*_\alpha x\neq f(x)$. Otherwise (if there is no such a $g_0$) choose any sequence $\langle O_n:n\in \omega \rangle$ of clopen subsets of $X_\alpha $ converging to $a^\alpha $.

Let ${\cal S}'={\cal S}''\cup \{\langle f(O_n):n\in \omega \rangle \}$. Look at $\{b^\alpha \}$ as a constant sequence and apply Lemma \ref{split} with ${\cal S}^{\prime }$ to get corresponding closed $A_0'$ and $A_1'$ subsets of $X_\alpha $ such that $(A_0',A_1')$ preserves each sequence in ${\cal S}'$. Note that by condition C1) we have that $X_\alpha \approx C$ so we can indeed use Lemma \ref{split}.

By construction there is a sequence $\{n_i:i\in \omega \}\subseteq \omega $ such that the subsequence $\langle f(O_{n_i}) :i\in \omega \rangle $ is entirely contained in $A_0'$ and converges to the point $b^\alpha$. Now we choose $b_i^\alpha \in f(O_{n_i})$ and let $a_i^\alpha=f^{-1}(b_i^\alpha)$ for $i\in \omega $.

There are two reasons for choosing the $b_i^\alpha $'s in this way. The first one is that $\{b_i^\alpha :i\in \omega \}$ still avoids ${\cal S}'$ and the second reason is that $\{ b_i^\alpha :i\in\omega \}$ also avoids the set $\{\langle g*_\alpha a_n^\alpha :n\in \omega \rangle:g\in G_\alpha \} $.

Thus if we let ${\cal S}={\cal S}'\cup \{\langle g*_\alpha a_n^\alpha :n\in \omega \rangle :g\in G_\alpha \}$ we can use Lemma \ref{split} again to get $A_0,A_1\subseteq X_\alpha $ closed such that $A_0\cup A_1=X_\alpha $, $A_0\cap A_1=\{b^\alpha \} $, $\{ b_i^\alpha :i\in \omega \}\subseteq A_0$ and $(A_0,A_1)$ preserves each sequence in ${\cal S}$.

If the hypothesis of C6) are satisfied (and hence the ones for C7) are not) and $x$ is the limit point of $\{z_n^\alpha :n\in \omega \}$ then we just apply Lemma \ref{split} directly to the set ${\cal S}''$ to get $A_0,A_1$ with $A_0\cap A_1=\{x\}$. If neither the hypothesis of C7) nor the ones for C6) are satisfied we do the same with an arbitrary $x\in X_\alpha$.

Identifying $X_\alpha$ with the cantor set $C$, we can now construct $Y(A_0,A_1,G) \subseteq X_\alpha \times 2^G$ as in the previous section and then identifying $2^G$ with $C$ we get $X_{\alpha +1}\subseteq X_\alpha \times C$.

It is clear from the construction that all the requirements are preserved. For C5) just use Lemma \ref{preserve}. Note that $X_\alpha$ has no isolated points by condition C8).

\section{Applications}\label{apps}

We assume $\Diamond$ throughout this section.

\begin{theorem}
There is a rigid Fedorchuck space.
\end{theorem}

\begin{proof}
In the main construction just let $G_\alpha =1$ for all $\alpha\leq \omega _1$.
\end{proof}

\begin{theorem}
There is an $h$-homogeneous Fedorchuck space.
\end{theorem}

\begin{proof} Let $\{(K_\gamma ,L_\gamma):\gamma \in\omega _1\} $ be an enumeration of all pairs of disjoint non-empty clopen subsets of $C\times C^{\omega _1}$. Let $G_0=1$ and assume $G_\alpha
'$ has been defined. Let $\gamma _0$ be the smallest $\gamma $ not considered at a previous stage for which both $K_\gamma $ and $L_\gamma $ have their supports contained in $\alpha $. Let $K=X_\alpha \cap \pi _\alpha K_{\gamma _0}$ and $L=X_\alpha \cap \pi _\alpha L_{\gamma _0}$. If both $K$ and $L$ are non-empty (and hence clopen in $X_\alpha $), apply Lemma \ref{free} to get $G_\alpha =\langle G_\alpha ',g\rangle $ acting freely on $X_\alpha $ and having $gK=L$. This clearly ensures that $X_{\omega _1}$ is $h$-homogeneous.
\end{proof}

\begin{theorem}
There is a zero-dimensional, hereditarily separable compact space $X$ which is $B$-homogeneous but not $h$-homogeneous.
\end{theorem}

\begin{proof}
Start with $G_0$ acting freely on $2^\omega $ such that any two disjoint clopen subsets are $G$-equivalent (for example one could use Lemma \ref{free} $\aleph _0$ times). Then keep $G_\alpha =G_0$ for all $\alpha \leq\omega _1$, in particular $G=G_0$. By Lemma \ref{nonequiexists}, $X_1$ has already two disjoint non-$G$-equivalent clopen subsets and by Lemma \ref{nonequipres} these remain non-$G$-equivalent (and hence non homeomorphic to each other) in $X$. However by Lemma \ref{base}, $X$ is still $B$-homogeneous.
\end{proof}

Under $MA+\neg CH$ any hereditarily separable compact space is first countable, so if it is also $B$-homogeneous it would have to be $h$-homogeneous by \cite{mat}. Thus our space could not have been constructed in $ZFC$. We don't know if there is a counterexample to Matveev's question in $ZFC$ (but again, if there is one, it can't be hereditarily separable).

\end{document}